\documentclass[12pt]{article}
\usepackage{latexsym,amssymb,amsfonts}
\usepackage{amsmath}
\usepackage{amsthm}
\usepackage{amscd}
\usepackage{graphicx}
\usepackage[arc,all]{xy}
\usepackage{tikz-cd}
\usepackage{braids}
\usepackage{hyperref}
\usepackage{pict2e}
\usepackage{ifpdf}

\theoremstyle{definition}
\theoremstyle{example}
\theoremstyle{remark}
\newtheorem{theorem}{Theorem}[section]
\newtheorem{proposition}{Proposition}[section]
\newtheorem{lemma}{Lemma}[section]
\newtheorem{definition}{Definition}[section]
\newtheorem{example}{Example}[section]
\newtheorem{remark}{Remark}[section]
\newtheorem{corollary}{Corollary}[section]

\begin{document}

\markboth{Uwe Kaiser}
{Jones and Kauffman bracket skein modules}

\title{Generalizing the relation between the Kauffman bracket and Jones polynomial \\
\
\\
\textsl{\small{Dedicated to the memory of Vaughn Jones}}}

\maketitle

\centerline{\author{UWE KAISER}}

\centerline{\textit{Department of Mathematics}}
\centerline{\textit{Boise State University, 1910 University Drive}}
\centerline{\textit{ Boise, ID 83725-1555, ukaiser@boisestate.edu
}}

\begin{abstract}
We generalize Kauffman's famous formula defining the Jones polynomial of an oriented link in $3$-space from his bracket and the writhe of an oriented diagram \cite{Ka}. Our generalization is an epimorphism between skein modules of tangles in compact connected oriented $3$-manifolds with markings in the boundary. Besides the usual Jones polynomial of oriented tangles we will consider graded quotients of the bracket skein module and Przytycki's $q$-analogue of the first homology group of a $3$-manifold \cite{P3}. In certain cases, e.\ g.\ for links in submanifolds of rational homology $3$-spheres, we will be able to define an epimorphism from the Jones module onto the Kauffman bracket module. For the general case we define suitably graded quotients of the bracket module, which are graded by homology. The kernels define new skein modules measuring the difference between Jones and bracket skein modules. We also discuss gluing in this setting. 

\end{abstract}

\centerline{\textit{Keywords:}{\ skein module, homology, 3-manifold}}

\centerline{Mathematics Subject Classification 2010: 57M25, 57M35, 57R42}

\section{Introduction}
The Kauffman bracket of an unoriented link projection $D$ is the element of $\mathbb{Z}[A^{-1},A]$ determined by
$$\langle D \rangle  = A \langle D_{cc}\rangle +A^{-1}\langle D_c\rangle, \quad \langle K\sqcup U\rangle = -A^2-A^{-2}, \quad \langle U\rangle =1$$
where $K_{cc}$, respectively $K_c$, is the counter-clockwise respectively clockwise smoothing at an unoriented crossing, and $U$ is the crossing-less projection of an unknot.  The notions of counter-clockwise and clockwise refer to opening the region swept by moving the overcrossing arc into the undercrossing arc at the crossing. 

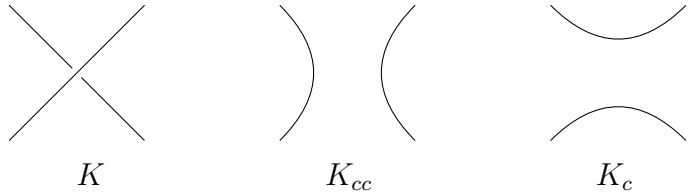
\setlength{\unitlength}{0.6cm}
\begin{picture}(5,5)
\put(2,0){\line(1,1){3}}
\put(5,0){\line(-1,1){1.4}}
\put(3.4,1.6){\line(-1,1){1.4}}
\put(3.5,-1){$K$}
\qbezier(8,0)(9.5,1.5)(8,3)
\qbezier(11,0)(9.5,1.5)(11,3)
\put(9,-1){$K_{cc}$}
\qbezier(14,0)(15.5,1.5)(17,0)
\qbezier(14,3)(15.5,1.5)(17,3)
\put(15,-1){$K_c$}
\end{picture}

\vskip 0.5in

The Jones polynomial of an oriented link is the element in $\mathbb{Z}[t^{-1/2},t^{1/2}]$ determined by the skein relations 
$$t^{-1}V(K_+)-tV(K_-)=(t^{1/2}-t^{-1/2})V(K_0), \quad V(K\sqcup U)=1$$
where $K_{\pm}$ and $K_0$ are defined by 

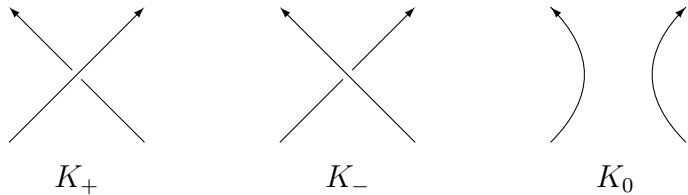
\setlength{\unitlength}{0.6cm}
\begin{picture}(5,5)
\put(2,0){\vector(1,1){3}}
\put(5,0){\line(-1,1){1.4}}
\put(3.4,1.6){\vector(-1,1){1.4}}
\put(3,-1){$K_+$}
\put(8,0){\line(1,1){1.4}}
\put(9.6,1.6){\vector(1,1){1.4}}
\put(11,0){\vector(-1,1){3}}
\put(9,-1){$K_-$}
\qbezier(14,0)(15.5,1.5)(14,3)
\qbezier(17,0)(15.5,1.5)(17,3)
\put(14.1,2.9){\vector(-1,1){0.1}}
\put(16.9,2.9){\vector(1,1){0.1}}
\put(15,-1){$K_0$}
\end{picture}

\vskip 0.5in
and $U$ here denotes the oriented unknot. 

It is the classical result by Kauffman \cite{Ka} (see also \cite{L}, Definition 3.6) that if $D$ is an oriented link projection corresponding to the oriented link $K$ then 
$$V(K)=\left((-A)^{-3w(D)}\langle \overline{D}\rangle\right)_{t^{1/2}=A^{-2}}\in \mathbb{Z}[t^{-1/2},t^{1/2}],$$
where $\overline{D}$ is the underlying unoriented diagram and $w(D)$ is the writhe of the oriented diagram. 

In this note we will generalize Kauffman's result to skein modules of tangles in oriented $3$-manifolds. In Definition 2.5 we define a new skein module, the relative Jones-bracket skein module. We will also discuss the compatibility with gluing.

\section{Statement of the main results}

Let $M$ be a compact connected oriented $3$-manifold, possibly with nonempty boundary $\partial M$. A \textit{marking} of $M$ is a $0$-dimensional submanifold $P$ of $\partial M$ consisting of an even number of points equipped with a non-zero vector tangent to $\partial M$ at each $x\in P$.  The marking is \textit{oriented} if half of the points are positively oriented and half of the points are negatively oriented. Let $\mathcal{M}$, respectively $\mathcal{M}'$, denote the set of oriented markings, respectively unoriented markings, of $M$. We should in fact use notation $\mathcal{M}(M)$ but the manifold is usually clear. Note that as zero-dimensional abstract manifold, $P$ is oriented null-bordant in the oriented setting, and null-bordant in the unoriented one. We will use the notation $P$ indicating the marking, so including the framing vectors. 

Sometimes we write $\mathcal{M}^*$ to indicate that a statement applies to $P\in \mathcal{M}$ or $P\in \mathcal{M}'$. For $P\in \mathcal{M}$ we let $P'$ denote its image in $\mathcal{M}'$, which is defined by forgetting the orientation. Note that each unoriented marking $Q$ can be oriented such that there is $P$ with $P'=Q$, i.\ e.\ $\mathcal{M}\rightarrow \mathcal{M}'$ is onto.  In the following we usually abbreviate $(M,\emptyset )=M$. 

For $P\in \mathcal{M}^*$ let $\mathcal{L}(M,P)$ denote the set isotopy classes of framed (oriented for oriented markings) tangles $K$ in $M$ bounding $P$ (in an oriented way for oriented markings). This means that $K$ is a union of annuli and bands. We have annuli embedded in the interior $\textit{Int}(M)$. If $b\cong I\times I$ is a band then $b\cap \partial M=(\partial I \times I)\cap \partial M$ consist of two arcs, which are perpendicular to the framing vectors at the two points $x,y$ in the boundary of the core $c$ of the band $I\times \{1/2\}$ and the framing of the normal bundle of $b$ in $M$ along $\partial c$ extends smoothly over $c$. Isotopy is ambient isotopy in $M$ fixing $\partial M$. Note that a diffeomorphism 
$(M_1,P_1)\rightarrow (M_2,P_2)$ will induce a bijective map $\mathcal{L}(M_1,P_1)\rightarrow  \mathcal{L}(P_2,M_2)$. 
For $P\in \mathcal{M}$ there is the natural surjective forget map $\mathcal{L}(M,P)\rightarrow \mathcal{L}(M,P')$.
If $P=\emptyset $ we also include the \textit{empty link} $K=\emptyset $, and apply skein relations to the empty link. 

In the following let $R=\mathbb{Z}$ for $\mathcal{M}$, and $R=\mathbb{Z}_2$ for $\mathcal{M}'$. If no coefficients are used then we refer to $\mathbb{Z}$.

Let $\partial: H_1(M,P;R)\rightarrow H_0(P;R)$ be the boundary homomorphism and \\ $h_1(M,P;R):=\partial^{-1}[P]$, where $[P]$ is the $R$-fundamental class of the $0$-manifold $P\in \mathcal{M}^*$. There is the natural surjective map 
$$\mathcal{L}(M,P)\rightarrow h_1(M,P;R).$$
assigning to each tangle its underlying homology class. The homology class of the empty link will be $0$. 
For $\alpha \in h_1(M,P;R)$ let $\mathcal{L}^{\alpha }(M,P)$ denote the set of isotopy classes of tangles with underlying homology class $a$. Then 
$$\mathcal{L}(M,P)=\bigsqcup_{\alpha \in h_1(M;R)} \mathcal{L}^{\alpha }(M,P).$$
Note that the sets $h_1(M,P;R)$ in general are not subgroups or $\mathbb{Z}/2\mathbb{Z}$-subspaces.

\remark It follows from the exact sequence of the pair $(M,P)$,
$$0\cong H_1(P;R)\rightarrow H_1(M;R)\rightarrow H_1(M,P;R)\rightarrow H_0(P;R)\rightarrow H_0(M;R)\cong R,$$
that $H_1(M)$ acts transitively and faithfully on $h_1(M;P)$, with $[P]$ mapping to $0$ under $H_0(P;R)\rightarrow H_0(M;R)$. In particular $|h_1(M,P)|=|H_1(M)|$ for every $P$. (Here, $|S|$ is the cardinality of a set $S$). In general, if we pick $\alpha \in h_1(M,P)$ then any $\beta \in h_1(M,P)$ has the form $\beta =\alpha + \gamma $ for a unique $\gamma \in H_1(M)$. Thus, if we are representing elements of $h_1(M,P)$ by tangles, we can choose a fixed tangle (e.g. without any closed components) and then represent each homology class just by adding closed components. 

\example Let $P\in \mathcal{M}^*$ be a marking of the compact $3$-ball $D^3$. Then $|h_1(M,P)|=1$. But if $H_1(M;R)$ is infinite then $|h_1(M,P)|$ is infinite. Exceptions for example are Lens-spaces or other rational homology $3$-spheres, possibly with $3$-balls or tubular neighborhoods of links removed but arbitrary markings allowed. On the other hand $|h_1(M,P)|$ is finite if $P\in \mathcal{M}'$.  

\vskip 0.1in

\definition Let $P\in \mathcal{M}$. A (homological) \textit{writhe invariant} is, for each $\alpha \in h_1(M,P;\mathbb{Z})$, a nonnegative integer $\omega (\alpha )\in \mathbb{Z}$ and a function
$$w=w_{\alpha }: \mathcal{L}^{\alpha }(M,P)\rightarrow \mathbb{Z}/2\omega (a)\mathbb{Z}$$
such that 
$$w(K_{\pm })=w(K_0)\pm 1,\quad w(K\amalg U)=w(K)$$
where $K_{\pm},K_0$ are defined locally in oriented $3$-balls in $M$ just like above with blackboard framings, and $U$ is a trivially framed unknot in a $3$-ball separated from $K$. 

\vskip 0.1in

\definition We say that the writhe function $w$ for $(M,P)$ has no writhe-indeterminacy if $\omega (\alpha )=0$ for all $\alpha \in h_1(M,P)$. We also say that $(M,P)$ has no writhe-indeterminacy if there exists a writhe function without indeterminacy.

\vskip 0.1in

\remark In considering skein related tangles some care is necessary. In fact, if we consider a \textit{skein ball} containing the framed tangle $K$ or $K_{\pm }$ then the normal vectors to the two disjoint bands in the $3$-ball both have to point upwards or both downwards. Otherwise a smoothing could obviously give rise to M\"obius bands instead of annuli. 

\vskip 0.1in

The inclusion $(M,P)\subset (M,\partial M)$ and the intersection pairing in $M$ define the bilinear pairing:
$$\iota : h_1(M,P)\otimes H_2(M)\rightarrow H_1(M,\partial M)\otimes H_2(M)\rightarrow \mathbb{Z}$$

\theorem For each $3$-manifold $M$ and $P\in \mathcal{M}$ choose a representative tangle $K^{\alpha }$ with homology class $\alpha \in h_1(M,P)$. Then there is defined the writhe map
$$\mathcal{L}(M,P)\rightarrow \mathbb{Z}/\omega (\alpha ),$$
where $\omega (\alpha )$ is the non-negative generator of the subgroup of $\mathbb{Z}$ defined by 
$$\iota (\alpha \otimes H_2(M))\subset \mathbb{Z}.$$

\vskip 0.1in

Theorem 2.1 is a generalization to tangles of Theorem 2.3 in \cite{P2}. We will give an independent  proof in section 3 based on Chapter 2 of \cite{K2}. 

\remark If $P=\emptyset $ then $M$ has no writhe-indeterminacy $\Longleftrightarrow M$ is a submanifold of a rational homology $3$-sphere $\Longleftrightarrow b_1(\partial M)=2b_1(M)$ for the first Betti-numbers $b_1$, see \cite{K1}, A.4.
In general, if $H_2(M)=0$ then there is no writhe-indeterminacy for all $P$.
Recall that the intersection pairing $\mathfrak{i}: \mathcal{F}H_1(M,\partial M)\otimes H_2(M)\rightarrow \mathbb{Z}$ is non-singular, where $\mathcal{F}H_1(M,\partial M)$ is the quotient by the torsion subgroup, and $H_2(M)$ is free abelian for each $3$-manifold (by Poincare duality and universal coefficient theorem for cohomology).
Therefore, if $0\neq \gamma \in H_2(M)$ let $\beta \in H_1(M,\partial M)$ such that $\mathfrak{i}(\beta \otimes \gamma)\neq 0$. Because
$$H_1(M,\partial M)\cong H^2(M)\cong [M,\mathbb{C}P^N]$$
for $N$ large and by Pontryagin-Thom construction we can realize $\beta $ as a $1$-dimensional oriented proper submanifold $C$ of $(M,\partial M)$. Define $P=\partial C$. Then $\beta $ is the image of an element $\alpha \in h_1(M,P)$ with $\iota (\alpha \otimes \gamma )\neq 0$. This shows $\omega (\alpha )\neq 0$. This shows that in fact for $H_2(M)\neq 0$ there exists $P$ such that $(M,P)$ has non-=trivial indeterminacy for some $\alpha $. In fact, we can find for each basis element of $H_2(M)$ a marking $P$ such that for some $\alpha \in h_1(M,P)$ we have $\omega (\alpha )=1$. 

On the other hand, if $\alpha $ is torsion then $\omega (\alpha )=0$ and the writhe on tangles in homology class $\alpha $ is defined in $\mathbb{Z}$.

\example Let $M=\Sigma \times I$ for $\Sigma $ a closed oriented surface. Then $H_2(M)\cong \mathbb{Z}$ generated by the surface $\Sigma \times \{1/2\}$. If $P=\{(x,0),(x,1)\}$ for some $x\in \Sigma $ then with $\alpha $ the homology class of $\{x\}\times I$ we get $\omega (\alpha )=1$ and the writhe function is trivial on links in this homology class. On the other hand, if $P=\{(x,0),(y,0)\}$ for $x\neq y$ in $\Sigma $ then $\iota $ is trivial on $h_1(\Sigma \times I,P)\otimes H_2(M)$. It follows that the writhe function has no indeterminacy in this case. 

\vskip 0.1in

Next we define two skein modules. 

\definition For $P\in \mathcal{M}$ let the \textit{Jones module} $\mathcal{J}(M,P)$ be the quotient of the free $\mathbb{Z}[A^{-1},A]$-module with basis $\mathcal{L}(M,P)$ by the submodule generated by framed Jones skein elements
$$A^4K_+-A^{-4}K_--(A^{-2}-A^2)K_0,$$
and 
$$K\sqcup U+(A^{-2}+A^2)K$$
and 
$$K^{(+)}-K,$$
where $U$ is a trivially framed unknot in a $3$-ball in $\textrm{Int}M$ separated from the tangle $K$, and $K^{(\pm )}$ will denote the tangle $K$ with the framing changed by one positive twist in any of the components. Note that this abuse of notation is justified by the framing relation. 

\vskip 0.1in

For $P\in \mathcal{M}$ let $P^o$ be the result of forgetting the framing of the marking. Let $\mathcal{L}^o(M,Q)$ denote the set of isotopy classes of oriented tangles bounding the oriented unframed marking $Q$. Then, for each framed marking $P$ such that $P^o=Q$, the obvious forget map $\mathcal{L}(M,P)\rightarrow \mathcal{L}^o(M,Q)$ is surjective. Now let $\mathcal{J}^o(M,Q)$ denote the quotient of the free $\mathbb{Z}[A^{-1},A]$-module by the submodule generated by the skein relations above for oriented tangles without the framing relation. 

\proposition The forget map $\mathcal{L}(M,P)\rightarrow \mathcal{L}^o(M,Q)$ defines an isomorphism of skein modules 
$$\mathcal{J}(M,P)\rightarrow \mathcal{J}^o(M,Q)$$
for each $P\in \mathcal{M}$ with $P^o=Q$. 
\proof This is immediate from the definitions and $K^{(+)}=K$ for the images of framed tangles in $\mathcal{L}(M,P)$. To define the inverse homomorphism, choose for each oriented link in $\mathcal{L}(M,P)$ an inverse in $\mathcal{L}(M,Q)$ and note that as an element of $\mathcal{J}$ the choice does not depend on the framing. $\blacksquare$
 
\vskip 0.1in
 
\remark (a) Note that we cannot add an empty object in a skein module of tangles if $P\neq \emptyset $ in a similar way as it is often done for skein modules of links. This is because in general there is not an obvious trivial tangle.  
 
\noindent (b) In the following we will not distinguish between the isomorphic modules $\mathcal{J}(M,P)$ and $\mathcal{J}^o(M,Q)$ and assume that we have a fixed choice of $P$ with $P^o=Q$. 

\vskip 0.1in

Similarly to the Jones module of tangles above we now define the \textit{(Kauffman) bracket module} $\mathcal{K}(M,Q)$ for $Q\in \mathcal{M}'$. This will be the quotient of the free $\mathbb{Z}[A^{-1},A]$-module with basis the set $\mathcal{L}(M,Q)$ of framed unoriented links in $M$ bounding $Q$ by the submodule generated by bracket skein relations for all tangles $K$ bounding $Q$:
$$K-AK_{cc}-A^{-1}K_c, \quad K\sqcup U+(A^{-2}+A^2)K$$

\definition For each $\mathbb{Z}[A^{-1},A]$-module $M$ and non-negative integer $n$ define:
$$M_n:=M/((A^{-n}-A^n)M)$$
be the $n$-\textit{reduced module}. 

\vskip 0.1in

Note that if $m\vert n$ then there is the projection epimorphism  $M_n\rightarrow M_m$.

\vskip 0.1in

\theorem Let $w$ be a writhe function for $(M,P)$. Then the associated Kauffman map, for $\alpha \in h_1(M,P)$:
$$\mathcal{L}^{\alpha }(M,P)\ni K\mapsto (-A)^{-3w(K)}K'\in \mathcal{L}(M,P')$$
defines an epimorphism of $\mathbb{Z}[A^{-1},A]$-modules:
$$\kappa : \mathcal{J}^{\alpha }(M,Q)\cong \mathcal{J}^{\alpha }(M,P)\rightarrow \mathcal{K}^{[\alpha ]_2}_{3\omega (\alpha )}(M,P')$$
where $P^o=Q$ and $P'\in \mathcal{M}'$ is the framed unoriented marking underlying $P$, $[\alpha ]_2\in h_1(M,P';\mathbb{Z}_2)$ is the $\mathbb{Z}_2$-reduction, and $\mathcal{J}^{\alpha }$ respectively $\mathcal{K}^{[\alpha ]_2}$ denote the submodules generated by tangles with the corresponding homology class.  

\remark (a) Note that the reduction of the Kauffman bracket module on the right hand side depends on $\alpha \in h_1(M,P)$ and not only on its $\mathbb{Z}_2$-reduction. So, in general, the reduction on the right hand side depends on the choice of lift of $[\alpha ]_2$ to $\alpha $. 

\noindent (b) Note that the skein modules of the compact $3$-ball $\mathcal{J}(D^3)$ and $\mathcal{K}(D^3)$ are each free $\mathbb{Z}[A^{-1},A]$-modules generated by the trivially framed unknot (respectively oriented unknot for $J^o(D^3,\emptyset )$. If we let $V$ denote the Jones polynomial of an oriented link and let $\langle \quad \rangle $ denote the bracket polynomial of a diagram (defining a framed link in $D^3$ using the blackboard framing) then we have the commutative diagram corresponding to Kauffman's formula:

\begin{center}
\begin{tikzcd}
\mathcal{J}^o(D^3)\cong \mathcal{J}(D^3) \arrow[r,"\kappa"] \arrow[rd,"V"] & \mathcal{K}(D^3) \arrow[d,"\langle \quad \rangle "]\\
                                                  &\mathbb{Z}[A^{-1},A]
\end{tikzcd}
\end{center}

\noindent (c) There is the action of the skein algebras of the $3$-ball, which in this case are identified with $\mathbb{Z}[A^{-1},A]$. In order to have consistency we include the empty link in the set of links and correspondingly adjust the skein relations, in particular having $\emptyset \sqcup U=-A^{-2}-A^2$ with the empty link corresponding to $1$. 
In this way the Kauffman map $\kappa $ can also be interpreted as a homomorphism of modules over the module of the $3$-ball. Note that $\kappa (1)=1$ using the corresponding identifications. 

\noindent (d) We can reduce from the formal case $\mathbb{Z}[A^{-1},A]$ by specifying $A$ to roots of unities, corresponding the values of the writhe indeterminacies, and we define epimomorphisms of complex vector spaces.   

\noindent (e) For $\alpha \in h_1(M,P)$ let $\mathfrak{k}^{\alpha }(M,P)$ denote the kernel of the corresponding Kauffman epimomorphism $\kappa ^{\alpha }$ as defined above. The study of this submodule of $\mathcal{J}^{\alpha }(M,P)$ seems to be of interest, in particular for $M=\Sigma \times I$ and $\Sigma $ an oriented surface, $P=\emptyset $, in which case there is no writhe indeterminancy, and the module is in fact an ideal of the skein algebra. 

\definition Let $\mathfrak{k}(M,P):=\bigoplus_{\alpha \in h_1(M,P)}\mathfrak{k}^{\alpha }(M,P)\subset \mathcal{J}(M,P)$ be the \textit{relative Jones-bracket module}.

\vskip 0.1in

\remark In the cases $(D^3,\emptyset )$ or $M=S^3$ the Kauffman epimorphism is an isomorphism. In general this is not the case: If $H_1(M)$ is infinite then $\mathfrak{k}(M,P)\neq 0$ just by grading reasons. Note that $\mathcal{J}^{\alpha }(M,P)$ and $\mathcal{K}^{\alpha }(M,P)$ are non-zero modules for each $\alpha \in h_1(M;R)$.

\example Consider $M=A\times I$ and $P=\emptyset $. In this case the skein modules are well-known and free by results of Turaev \cite{T} and Przytycki \cite{P4}. Note that $H_1(M)\cong \mathbb{Z}$ and $H_1(M,\mathbb{Z}_2)\cong \mathbb{Z}/2\mathbb{Z}$. Consider $a=0\in H_1(M)$. Then a basis of this part of the skein module is determined by sequences $(a_1,\ldots ,a_n)$ with $\sum_{i=1}^na_i=0$ and $a_1\geq a_2\geq \cdots \geq a_n$ corresponding to the components of a link with descending components and separated by the radius function of the annulus. The bracket module is isomorphic to $\mathbb{Z}[A^{-1},A][z]$ with $z$ being represented by the core circle of the annulus. Now consider basis elements of the Jones module corresponding to $(1,-1)$ and $(1,1)$. Both of these map to the same basis element of bracket module corresponding to $z^2$.  

\vskip 0.1in

\corollary Suppose that $\omega (\alpha )=0$ for all $\alpha \in h_1(M,P)$. Then the above epimorphisms combine to the epimorphism:
$$\kappa: \mathcal{J}(M,Q)\rightarrow \mathcal{K}(M,P')$$
and kernel $\mathfrak{k}(M,P)$. 

\remark The assumption of the corollary is true for example when $P=\emptyset $ and $M$ is a submanifold of a rational homology $3$-sphere, or if $M=\Sigma \times I$ for an oriented surface $\Sigma $
with $P\subset \partial \Sigma \times I$. 
  
\vskip 0.1in

\section{Proof of Theorem 2.1}
Even though not necessary it is helpful to restate Theorem 2.1 in the language of skein modules. 
We will give a different proof though using ideas from Chapter 3 of \cite{K1}. Just like in \cite{K1}, oriented bordism of tangles in $M\times I$ corresponds to homology in the group $h_1(M,P)$.

\definition (see Definition 2.1 in \cite{P3}) Given $(M,P)$ as above the \textit{Przytycki skein module} $\mathcal{P}(M,P)$ is the quotient of the free 
$\mathbb{Z}[q^{-1},q]$-module with basis $\mathcal{L}(M,P)$ by the submodule generated by the relations
$$K_+-qK_0,\quad K^{(+)}-qK$$

\remark (a) For $P=\emptyset $, the skein module $\mathcal{P}(M)$ is denoted $\mathcal{S}_2(M;\mathbb{Z}[q^{\pm 1}],q)$ in \cite{P3}.

\noindent (b) We say that the writhe invariant is universally defined for $(M,P)$ if $\omega (\alpha )=0$ for all $\alpha \in h_1(M,P)$.

\noindent (c) The Przytycki skein module has been called the $q$-analogue of the first homology group by Przytycki in \cite{P3} and computed by him in the link case. In \cite{Ka2} the $q$-analogue of the fundamental group of an oriented $3$-manifold has been computed.

\vskip 0.1in

The following equivalent of Theorem 2.1 generalizes Przytycki's Theorem 2.6 \cite{P3}. 

\theorem The choice of framed links $K_{\alpha }$ for each $\alpha \in h_1(M,P)$ defines the isomorphism 
$$\mathcal{P}(M,P)\rightarrow \bigoplus_{\alpha \in h_1(M,P)}\mathbb{Z}/(2\omega (\alpha )\mathbb{Z})$$

\proof By applying the first relation to the positive crossing in $K^{(+)}$
we get in $\mathcal{P}(M,P)$, $K^{(+)}=qK\sqcup U=qK$, which implies that $K\sqcup U=K$ holds in the skein module. 

\vskip -0.1in

\setlength{\unitlength}{0.6cm}
\begin{picture}(5,5)
\put(8,0){\line(1,1){2.0}}
\put(10,1){\line(-1,1){0.4}}
\put(9.4,1.6){\vector(-1,1){1.4}}
\qbezier(10,2.0)(11.5,4.1)(13,2.0)
\qbezier(10,1.0)(11.5,-0.5)(13,1.2)
\qbezier(13,2.0)(13.22,1.6)(13,1.2)
\put(9,-1.5){$K^{(+)}$}
\end{picture}

\vskip 0.5in

We also see from the two relations that $K_+=K_0^{(+)}$, similarly $K_-=K_0^{(-)}$. (This is exactly the spin-statistics observation made by Kauffman in \cite{Ka2}, p.\ 12).  It follows that we can arbitrarily change a framed tangle in its oriented bordism class as long as we keep track of framing changes corresponding to band attachments. Thus fixing a framed link $K_{\alpha }$ for given $\alpha \in h_1(M,P)$, a tangle $K$ in $\mathcal{P}(M,P)$ 
with homology class $\alpha $ is determined by comparing the resulting framing after an oriented bordism of $K$ in $M\times I$ with the framing of $K_{\alpha }$. In fact, just like observed in Lemma 2.2.4 of \cite{K1} the bordism group of the framed $3$-ball acts transitively on framed bordism in $M\times I$. Note that we can add a trivially framed circle in a separated $3$-ball to $K_{\alpha }$, which will not change the image in the skein module, and then change framings of this unknot to change the framings of $K_{\alpha }$. The indeterminancy of the power of $q$ resulting from this is determined by the normal Euler class of the surface in $M\times S^1$ resulting from a self-bordism of $K_{\alpha }$ into $K_{\alpha }\sqcup U$, with $U$ capped off by a disk, see Lemma 2.4.5 in \cite{K1}. This normal Euler class is the self-intersction number of the bordism surface with itself. Now note that the self-bordism represents an element in the inverse image under the boundary operator 
$H_2(M\times S^1,P\times S^1)\rightarrow H_1(P\times S^1)$ of $[P]\times [S^1]$, where as before $[\quad ]$ denotes the fundamental class. Note that the K\"unneth formula applied to the pair $(M,P)\times (S^1,\emptyset )$ gives:
$$H_2(M\times S^1,P\times S^1)\cong H_2(M,P)\oplus H_1(M,P).$$
From the long exact homology sequence of $(M,P)$ we get $H_2(M,P)\cong H_2(M)$ and the condition about the inverse image above implies that the surface represents an element in 
$$H_2(M)\oplus h_1(M,P).$$
The self-intersection number of this surface, therefore the possible change in framing, is given by twice the intersection pairing applied to the element of $H_2(M)$ and $h_1(M,P)$ representing the self-bordism. Note that surjectivity follows from 
transitivity of the action of the framed skein module of the $3$-ball. $\blacksquare$

\vskip 0.1in

Of course if $\alpha \in h_1(M,P)$ is torsion then $\omega (\alpha )=0$. Therefore, if $P=\emptyset $ and $M$ is a submanifold of a rational homology $3$-sphere then $\omega (\alpha )=0$ for all $\alpha $, in this case $H_1(M)=h_1(M)$ is torsion. See \cite{K1} for more details about this class of manifolds called Betti-trivial in \cite{K1}. 

\vskip 0.1in

\example Let $M=\Sigma \times I$ for $\Sigma $ an oriented compact surface. This is a submanifold of a rational homology $3$-sphere so the writhe invariant is universally defined for $P=\emptyset $. Note that intersections of $1$-cycles with $2$-cycles in the submanifold cannot be non-trivial because we would get the same result in a rational homology sphere, into which it embeds. The situation is more subtle when $P\neq \emptyset $. If $\Sigma $ is a closed surface and for example 
$P=\{(x,0),(x,1)\}$ for some $x\in \Sigma $ then $\omega (\alpha )=1$ because the surface $\Sigma \times \{1/2\}$ will have intersection number $\pm 1$ with e.g. the class in $h_1(\Sigma \times I,P)$ defined by $\{x\}\times I$. 
But if $P\subset \Sigma \times \{i\}$ for $i=0$ or $i=0$ then intersection numbers will be $0$ and the writhe is universally defined. If $\partial \Sigma \neq \emptyset $ then $P$ can be isotoped into $\partial \Sigma \times I$, and at the same time $H_2(\Sigma \times I)=0$ so our intersection pairing is trivial. Note that for $\Sigma =S^2$ we get the familiar consequence of the light bulb theorem. In this case for a closed framed loop in $S^2\times S^1$ or equivalently a tangle in $S^2\times I$ intersecting $S^2\times \{1/2\}$ in a single point  the framing can be changed by sign through isotopy across the $2$-sphere, see \cite{C1}. Our case generalizes this observation to an oriented surface, which represents an oriented bordism changing the framing correspondingly. 

\vskip 0.1in

\section{Proof of Theorem 2.2}

First not that it is well-known that the bracket relations imply that $[L^{(+)}]=-A^3[L]=(-A)^3[L]$ for any unoriented tangle in $(M,P)$ and $[\quad ]$ denotes the class in the bracket module. Therefore
$$\kappa (K^{(+)})=\kappa (K)\in \mathcal{K}(P,M).$$
This follows using $w(K^{(+)})=w(K)+1$:
$$\kappa (K^{(+)})=(-A)^{-3w(K^{(+)})} [{K^{(+)}}']=(-A)^{-3w(K)-3}(-A)^3[K']=\kappa (K),$$
where in general $L'$ is the unoriented framed tangle determined by the oriented framed tangle $L$. 
The writhe invariant suffices to pick from the orientation, in combination with the framing, a framing of the unoriented link, which is independent of the framing of the original tangle. Framings can be understood as \textit{local} projections, so not surprisingly there are obstructions to be able to combine the local information to a global map like $\kappa $.

\vskip 0.1in

To finish the proof of well-definedness of $\kappa $ we calculate $\kappa $ on $A^4K_+-A^{-4}K_-$ and compare with 
the mapping applied to $(A^{-2}-A^2)K_0$. We find first using the property of the writhe and the bracket relations, and 
noting that if $K_+'=:K$ and $K_-'=:L$ then $K_{cc}=L_c=K_0'$ and $K_c=L_{cc}$, the following holds:
\begin{align*}
\kappa (A^4K_+-A^{-4}K_-)&=A^4(-A)^{-3w(K_+)}[K_+']-A^{-4}(-A)^{-3w(K_-)}[K_-']\\
&=A^4(-A)^{-3w(K_0)-3}[K_+']-A^{-4}(-A)^{-3w(K_0)+3}[K_-']\\
&=(-A)^{-3w(K_0)}\left(-A[K_+']+A^{-1}[K_-']\right)\\
&=(-A)^{-3w(K_0)}\left(-A\left(AK_{cc}-A^{-1}K_c\right)+A^{-1}\left(AK_c+A^{-1}K_{cc}\right)\right)\\
&=(-A)^{-3w(K_0)}(-A^2+A^{-2})K_{cc}\\
&=(A^{-2}-A^2)\kappa (K_0)
\end{align*}

In the above calculation we did not indicate that $w$ is not necessarily taking values in $\mathbb{Z}$. The calculations though are valid in the quotients of the Kauffman bracket depending on the integral homology class of $K_{\pm},K_0$, as defined above.  

The surjectivity of $\kappa $ follows immediately from the definitions and the surjectivity of the marking map $\mathcal{M}\rightarrow \mathcal{M}'$. $\blacksquare$

Note that by the application of bracket relations to a tangle will keep the mod $2$ homology class, but in terms of picking preimages in the Jones skein module the integral homology class of a preimage is not naturally defined.

\section{Compatibility with gluing}

We decided to prove Theorem 2.1 for tangles because our homomorphisms are compatible with gluings, even though the behavior of writhe indeterminacy under gluing adds technical difficulties. 

Let $(M_1,P_1)$ and $M_2,P_2)$ be two marked oriented $3$-manifolds. Let $\partial_0(M_i)$ be compact submanifolds of $\partial M_i$ for $i=1,2$. Let $h: \partial_0M_1\rightarrow \partial_0M_2$ be an orientation reversing diffeomorphism. Let 
$P_i^0:=P_i\cap \partial_0M_i\subset \textrm{int}\partial_0M_i$ for $i=1,2$. Let $h(P_1^0)=\overline{P_2^0}$, where $\overline{\ }$ denotes the opposite orientation. We also assume that the derivative of $P$ at each point $x\in P_1^0$ maps the marking vector at that point to the corresponding marking vector at $h(p)$. Let $P:=(P_1\setminus P_1^0)\sqcup (P_2\setminus P_2^0)\subset \partial M_1\cup_h \partial M_2=\partial M$ where $M$ is the boundary connected sum   
$M_1\cup_h M_2$. Note that there are well-defined maps of sets of framed tangles defined by gluing the tangles using $h$:
$$\mathcal{L}(M_1,P_1)\sqcup \mathcal{L}(M_2,P_2)\rightarrow \mathcal{L}(M,P).$$
Note that gluing using $h$ also defines the map
$$\mathfrak{h}: h_1(M_1,P_1)\times h_1(M_2,P_2)\rightarrow h_1(M,P).$$
Note that the $h_1(M_i,P_i)$ in general are only subsets of $H_1(M_i,P_i)$, $i=1,2$. The above map is the restriction of the corresponding homomorphism  
$$H_1(M_1,P_1)\oplus H_1(M_1,P_2)\rightarrow H_1(M,P_1\cup P_2)$$
in the relative Mayer-Vietoris sequence (see e.\ g.\ \cite{S}, p.\ 187), where we have identified $M_i\subset M$, and the image is in the subgroup $H_1(M,P)$ because the boundaries are glued. Note that $P=(P_1\cup P_2)\setminus (P_1\cap P_2)$ and there is the exact sequence of the triple $P\subset P_1\cup P_2\subset M$ (see \cite{S}, p.\ 201):
$$0=H_1(P_1\cup P_2,P) \rightarrow H_1(M,P)\rightarrow H_1(M,P_1\cup P_2)\rightarrow H_0(P_1\cup P_2,P).$$
Note that the same homological discussion as above works for coefficients in $\mathbb{Z}_2$. 

\lemma Let $\alpha_i \in h_1(M_i,P_i)$ for $i=1,2$. Then $\omega (\mathfrak{h}(\alpha_1,\alpha_2))$ divides the gcd of $\omega (\alpha_1)$ and $\omega (\alpha_2))$.

\proof Note that 
$$\iota (\mathfrak{h}(\alpha_1,\alpha_2)\otimes H_2(M_i))=\iota (\alpha_i \otimes H_2(M_i)),$$
where we have identified $H_2(M_i)$ with its image in $H_2(M)$ under the inclusions $M_i\subset M$. 
Therefore there are more intersections to be considered for $\mathfrak{h}(\alpha_1,\alpha_2)$ than coming from intersections with surfaces in homology classes in $H_2(M_i)$ for $i=1,2$. $\blacksquare$

\theorem Let $(M_i,P_i)$ for $i=1,2$ and $h$ be as above. Let $\alpha_i\in h_1(M_i,P_i)$ for $i=1,2$. Then, for given writhe functions for $(M_1,P_1), (M_2,P_2)$, we can choose a writhe function for $(M,P)$ such that there is the commutative diagram of skein modules with the natural gluing maps for skein modules horizontally:

\begin{center}
\begin{tikzcd}
\mathcal{J}^{\alpha_1}(M_1,P_1)\otimes \mathcal{J}^{\alpha_2}(M_2,P_2)\arrow[r,"\mathfrak{g}_{\mathcal{J}}"] \arrow[d,"\kappa^{\alpha_1}\otimes \kappa^{\alpha_2}"] & \mathcal{J}^{\mathfrak{h}(\alpha_1,\alpha_2)}(M,P)\arrow[d,"\kappa^{\mathfrak{h}(\alpha_1,\alpha_2)}"]\\
\mathcal{K}^{[\alpha_1]_2}_{3\omega (\alpha_1)}(M_1,P_1')\otimes \mathcal{K}^{[\alpha_2]_2}_{3\omega (\alpha_2)}(M_2,P_2') \arrow[r,"\mathfrak{g}_{\mathcal{K}}"] \arrow[d] & \mathcal{K}^{[\mathfrak{h}(\alpha_1,\alpha_2)]_2}_{3\omega (\mathfrak{h}(\alpha_1,\alpha_2))}(M,P')\\
\mathcal{K}^{[\mathfrak{h}(\alpha_1,\alpha_2)]_2}_{3\textrm{gcd}(\omega (\alpha_1),\omega (\alpha_2))}(M,P')\arrow[ru]
\end{tikzcd}
\end{center}

\proof The vertical left bottom arrow is defined by gluing and the fact that the $\textrm{gcd}$ divides $\omega (\alpha_i)$ for $i=1,2$, the arrow in the second row is defined by commutativity of the triangle. The subtle point is to choose the writhe function used to define $\kappa^{\mathfrak{h}(\alpha_1,\alpha_2)}$ from the representing framed tangle by gluing the representative tangles, which have been used to define the writhe maps $\kappa^{\alpha_i}$ for $(M_i,P_i)$. The result then follows from a gluing bordisms argument, which will imply $w(K_1\cup_hK_2)=w(K_1)+w(K_2)$, where $K_1\cup_hK_2$ is the result of gluing two framed oriented angles $K_1,K_2$ using $h$. Note that
\begin{align*}
\mathfrak{g}_{\mathcal{K}}(\kappa^{\alpha_1}\otimes \kappa^{\alpha_2})(K_1\otimes K_2)&=\mathfrak{g}_{\mathcal{K}}((-A)^{-3w(K_1)}[K_1']\otimes (-A)^{-3w(K_2)}[K_2'])\\
\quad &=(-A)^{-3(w(K_1)+w(K_2))}[K_1'\cup_hK_2']
\end{align*}
and correspondingly
\begin{align*}
\kappa^{\mathfrak{h}(\alpha_1,\alpha_2)}\mathfrak{g}_{\mathcal{J}}(K_1\otimes K_2)&=(-A)^{-3w(K_1\cup_hK_2)}
[(K_1\cup_hK_2)']\\
\quad &=(-A)^{-3w(K_1\cup_hK_2)}[(K_1'\cup_hK_2)']
\end{align*}
because gluing and forgetting orientations commute. In the above formulas $w$ denotes three different writhe functions depending on the choices of framed tangles representing $\alpha_i$ for $i=1,2$ and $\mathfrak{h}(\alpha_1,\alpha_2)$, which is really the gluing of homology classes. 
 $\blacksquare$

\vskip 0.1in

\remark The connection between the Jones and bracket module is Przytycki's module with 
$$\mathcal{P}^{\alpha }(M,P)\cong \mathbb{Z}/(2\omega (\alpha )\mathbb{Z}),$$
with the isomorphism defined by a choice of framed link in the homology class $\alpha \in h_1(M,P)$.
Gluing of framed tangles as defined above induces the homomorphism
$$\mathcal{P}^{\alpha_1}(M_1,P_1)\otimes \mathcal{P}^{\alpha_2}(M_2,P_2)\rightarrow \mathcal{P}^{\mathfrak{h}(\alpha_1,\alpha_2)}(M,P)$$ and induces with corresponding choices of representing framed tangles:

\begin{align*}
\mathbb{Z}/(2\omega (\alpha_1)\mathbb{Z})\otimes \mathbb{Z}/(2\omega (\alpha_2)\mathbb{Z})  & \rightarrow 
\mathbb{Z}/(2\textrm{gcd}(\omega (\alpha_1),\omega (\alpha_2)\mathbb{Z})\\
\quad & \rightarrow \mathbb{Z}/(2\omega (\mathfrak{h}(\alpha_1,\alpha_2))\mathbb{Z})
\end{align*}

(choosing the representative framed tangle for $\mathfrak{h}(\alpha_1,\alpha_2)$ by gluing the representative tangles for $\alpha_1$ and $\alpha_2$). Note that we do not claim that we can choose the writhe function for the glued homology class (equivalently the framed links in homology classes) for all choices of $\alpha_1,\alpha_2$ uniformly.

Recall that $M_1\cap M_2=\partial_0M$ after identifying $M_i\subset M$ for $i=1,2$.

\theorem Suppose that the boundary homomorphism $\partial: H_2(M)\rightarrow H_1(\partial_0M)$ in the exact sequence of $(M,\partial_0M)$ is trivial. Then, for given writhe functions for $(M_i,P_i)$ defining $\kappa^{\alpha_i}$ for $i=1,2$, there is defined a writhe function for $(M,P)$, $P\in \mathcal{M}$, defining $\kappa^{\alpha }$, such that the conclusion of Theorem 5.1 holds for all $\alpha_1, \alpha_2$ uniformly. 

\proof It follows from the assumption that the homomorphism $j$ in the Mayer-Vietoris sequence

\begin{center}
\begin{tikzcd}
H_2(M)\arrow[r,"\partial"] & H_1(\partial_0M)\arrow[r,"j"] & H_1(M_1)\oplus H_1(M_2)\arrow[r] & H_1(M)
\end{tikzcd}
\end{center}

is injective. Now $H_1(M_i)$ ($i=1,2$) and $H_1(M)$ act transitively and faithfully on the corresponding sets 
$h_1(M_i,P_i)$ ($i=1,2)$) and $h_1(M)$. Therefore the injectivity above implies that
$$\mathfrak{h}: h_1(M_1,P_1)\times h_1(M_2,P_2)\rightarrow h_1(M,P)$$
is injective. So for each $\alpha \in h_1(M,P)$ in the image of $\mathfrak{h}$ we can find unique $\alpha_i\in h_1(M_i,P_i)$ for $i=1,2$ such that $\mathfrak{h}(\alpha_1,\alpha_2)=\alpha $. Now for given framed tangles with homology classes $\alpha_i$ for $i=1,2$ we can construct the tangle in $(M,P)$ by gluing to represent the homology class $\alpha $. 
For all $\alpha $ not in the image of $\mathfrak{h}$ we choose the representing framed tangle arbitrarily. This defines a writhe function with the property claimed above. $\blacksquare$

\example A natural gluing example is the gluing of two handlebodies $\mathcal{H}_{g,i}$, $i=1,2$ along their  boundaries. We think of $\mathcal{H}_g=D_g \times I$, where $D_g$ is the $g$-holed disk and $\Sigma_g=\partial \mathcal{H}_g$. Note that for each $P$ in this case $h_1(\mathcal{H}_g),P)$ is bijective with $H_1(D_g)=\mathbb{Z}^g$. But $(\mathcal{H}_g,P)$ has no writhe-indeterminacy because $H_2(\mathcal{H}_g)=0$. Note that Mayver-Vietoris sequence extended to the left gives the following:

\begin{center}
\begin{tikzcd}
H_2(\mathcal{H}_{g,1})\oplus H_2(\mathcal{H}_{g,2})\arrow[r] & H_2(M)\arrow[r,"\partial"] & H_1(\Sigma_g),
\end{tikzcd}
\end{center}

which shows that $\partial $ is injective and therefore $\partial \neq 0$ for $H_2(M)\neq 0$. Thus if $H_2(M)=0$, Theorem 5.2 applies to a Heegaard splitting. 

\section{Final remarks}
Skein modules were introduced by Przytycki \cite{P1}, \cite{P2} for oriented $3$-manifolds and Turaev for links in surfaces \cite{T}. There has been much recent progress in the field, see e.\ g.\ the introduction of \cite{GJS}. Gilmer and Masbaum studied the bracket module of $\Sigma \times S^1$ for $\Sigma $ a closed oriented surface \cite{GM}. Frohman, Kania-Bartoszynska and Le developed the representation theory of the bracket module of $\Sigma \times I$ for $\Sigma $ an oriented surface. Gunningham, Jordan and Safranov proved a conjecture by Witten about finite-dimensionality of the bracket module of a closed oriented $3$-manifold over $\mathbb{C}(A)$ using the theory of factorization algebras.  The skein relations of a skein module are usually deduced from the $R$-matrix of a quantum group, and in the case of the quantum group of $\textrm{sl}(2)$ define the Jones relations, while approaches to the bracket using matrices can be found e.\ g.\ in \cite{Ka2}. In this article we defined a \textit{natural} homomorphism between the two modules, generalizing Kauffman's famous formula \cite{Ka}. In this way many results about the bracket module naturally give rise to results about Jones skein modules. The study of the kernels seems to be an interesting future project.

\end{document}